\documentclass{amsart}
\usepackage{amssymb,amsmath,amsfonts,amsthm,amscd,latexsym}
\usepackage{graphics}
\usepackage{psfrag}
\newcommand{\cal}[1]{\mathcal{#1}}
\theoremstyle{plain}
\newtheorem{lemma}{Lemma}[section]
\newtheorem{theorem}{Theorem}
\newtheorem*{cor}{Corollary}

\newtheorem{proposition}[lemma]{Proposition}
\newtheorem{corollary}[lemma]{Corollary}
\parskip=\bigskipamount

\let\egthree=\phi
\let\phi=\varphi
\let\varphi=\egthree

\usepackage[dvips]{graphicx} 
\usepackage[pdflinkmargin=5pt,pdfstartview={FitBH -32768},pdfpagemode=None,plainpages=true]{hyperref}
\newtheorem{lem}[lemma]{Lemma}
\newtheorem{kor}[lemma]{Corollary}
\newtheorem*{theo*}{Theorem}

\theoremstyle{definition}
\newtheorem{definition}[lemma]{Definition}

\def\IR{\mathbb{R}}

\def\IC{\mathbb{C}}

\begin{document}
\title[Topological properties of Reeb orbits]%
{Topological properties
of Reeb orbits on boundaries of star-shaped domains
in $\mathbb{R}^4$}
\author{Stefan Hainz, Ursula Hamenst\"adt}
%Mathematisches Institut der Universit\"at Bonn,\\
%Beringstra\ss{}e 1, D-53115 Bonn}
\thanks %{\it e-mail address:} ursula@math.uni-bonn.de\\
{Partially supported by DFG-SPP 1154.}

\date{February 7, 2010}

\begin{abstract} Let $B^4$ be the compact unit ball in 
$\mathbb{R}^4$ 
with boundary $S^3$.
Let $\gamma$ be a knot on $S^3$ which is transverse
to the standard contact structure.
We show that if there is an immersed symplectic disc
$f:(D,\partial D)\to (B^4,\gamma)$ then
${\rm lk}(\gamma)=2{\rm tan}(f)-1$ where ${\rm lk}(\gamma)$
is the self-linking number of $\gamma$ and ${\rm tan}(f)$ is the
tangential self-intersection number of $f$.
We also show that if $E\subset \mathbb{C}^2$ is 
compact and convex, with smooth boundary $\Sigma$, and if 
the principal curvatures of $\Sigma$ are suitably pinched then
the self-linking number of 
a periodic Reeb orbit on $\Sigma$ of Maslov index 3
equals $-1$.
\end{abstract}

\maketitle

\section{Introduction}

Consider the four-dimensional euclidean space
$\mathbb{R}^4$ with the standard \emph{symplectic form}
defined in standard coordinates by
$\omega_0=\sum_{i=1}^2dx_i\wedge dy_i$. This symplectic
form is the differential of the one-form
\[\lambda_0=\frac{1}{2}\sum_{i=1}^2(x_idy_i-y_idx_i).\]
For every bounded domain $\Omega\subset \mathbb{R}^4$
which is star-shaped
with respect to the origin $0\in
\mathbb{R}^4$, with smooth boundary $\Sigma$,
the restriction $\lambda$ of $\lambda_0$ to $\Sigma$
defines a smooth \emph{contact form} on $\Sigma$. This means that
$\lambda\wedge d\lambda$ is a volume form on $\Sigma$.

Let $\xi={\rm ker}(\lambda)$ be the contact bundle.
Each \emph{transverse knot} $\gamma$ on $\Sigma$, i.e. an 
embedded smooth closed curve
on $\Sigma$ which is everywhere transverse to $\xi$,
admits a canonical orientation determined by the requirement that
$\lambda(\gamma^\prime)>0$. To such an oriented transverse knot
$\gamma$  
we can associate
its \emph{self-linking number} ${\rm lk}(\gamma)$ which is defined
as follows. Let $S\subset \Sigma$ be a \emph{Seifert surface} for
$\gamma$, i.e. $S$ is a smooth embedded oriented 
surface in $\Sigma$ whose oriented 
boundary equals $\gamma$. Since $\gamma$ is transverse to $\xi$, 
there is a natural identification of 
the restriction to $\gamma$ of the oriented normal
bundle of $S$ in $\Sigma$ with
a real line subbundle $N_S$ of the contact bundle
$\xi\vert\gamma$. Then $N_S$ defines a
trivialization of the oriented two-plane bundle
$\xi\vert \gamma$. The self-linking number
${\rm lk}(\gamma)$ of $\gamma$ is the winding number 
with respect to $N_S$ of a
trivialization of $\xi$ over $\gamma$ which extends to a
trivialization of $\xi$ on $\Sigma$. 

Eliashberg \cite{eliash2}
showed that the self-linking number of a transverse
knot in $\Sigma$ is always an odd integer. If $g$ denotes 
the \emph{Seifert genus} of $\gamma$, i.e. the smallest genus of a
Seifert surface for $\gamma$, 
then we have
${\rm lk}(\gamma)\leq 2g-1$ 
(Theorem 4.1.1 of \cite{eliash2}). 
Eliashberg also constructed
for every $k\geq 1$ a transverse unknot $\zeta$ 
of self-linking number ${\rm lk}(\zeta)=-2k-1$
in the standard unit three-sphere $S^3\subset \mathbb{R}^4$.
This indicates that in general we can not expect additional
relations
between the self-linking number of a transverse knot $\gamma$
on $\Sigma$ and purely topological invariants of $\gamma$.

Our first goal is to relate the self-linking number of a 
(canonically oriented) 
transverse knot $\gamma$ on the unit three-sphere
$S^3\subset \mathbb{R}^4$ to
the symplectic topology of the closed 
unit ball $B^4\subset \mathbb{R}^4$.  For this let 
$D\subset \mathbb{C}$ be the closed unit disc with 
oriented boundary
$\partial D$ and let  
$f:(D,\partial D)\to (B^4,S^3)$ be a smooth  immersion with  
$f^{-1}(S^3)=\partial D$.
If all self-intersections of $f(D)$ are
transverse then 
the \emph{tangential index} ${\rm tan}(f)$ of $f$ is the
number of self-intersection points of $f$ 
counted with signs and multiplicities.
The disc $f:D\to B^4$ is called \emph{symplectic} 
if for every $x\in {D}$ the 
restriction of the symplectic form $\omega_0$ to $
df(T_x{D})$ does not
vanish and defines the usual orientation of $D$. We show

\begin{theorem}\label{thm1}
Let $\gamma$ be a transverse knot on the boundary $S^3$ of 
the compact unit ball $B^4\subset \mathbb{R}^4$.
If $\gamma$ bounds an immersed symplectic disc
$f:(D,\partial D)\to (B^4,\gamma)$
then ${\rm lk}(\gamma)=2 {\rm tan}(f)-1$.
\end{theorem}

Now let $\Sigma$ be the boundary 
of an arbitrary bounded domain $\Omega\subset \mathbb{R}^4$
which is star-shaped with respect to the origin, with smooth boundary.
The \emph{Reeb vector field} of the contact form $\lambda$ 
is the smooth vector field $X$ on $\Sigma$ defined by $\lambda(X)=1$
and $d\lambda(X,\cdot)=0$. 
Rabinowitz \cite{Ra78} (see also \cite{W79}) showed 
that the \emph{Reeb flow} on $\Sigma$
generated by the Reeb vector field $X$ admits periodic orbits. 
Dynamical properties of the Reeb flow on $\Sigma$ are
related to properties of $\Omega$ viewed as a symplectic manifold.
The proof of the following corollary is similar to the 
proof of Theorem \ref{thm1}.

\begin{cor}\label{reeb}
Let $\gamma$ be a periodic Reeb orbit on the boundary
$\Sigma$ of a star-shaped domain $\Omega\subset \mathbb{R}^4$
with compact closure $C$. If $\gamma$ bounds an immersed
symplectic disc $f:(D,\partial D)\to (C,\gamma)$ then
${\rm lk}(\gamma)=2{\rm tan}(f)-1$.
\end{cor}

Note that by a result of
Hofer, Wysocki and Zehnder \cite{HWZ96}, 
there is always 
a periodic Reeb orbit of self-linking number $-1$ on $\Sigma$ which
is unknotted.

Even though the radial diffeomorphism $\Psi:
S^3\to \Sigma$ maps the contact bundle of $S^3$ to 
the contact bundle of $\Sigma$ and hence maps a transverse
knot $\gamma$ on $S^3$ to a transverse knot $\Psi\gamma$ 
on $\Sigma$, the corollary
is not immediate from 
Theorem \ref{thm1}. Namely,  in general 
the radial diffeomorphism $\Psi$ does not extend to a 
symplectomorphism $B^4\to C$ 
and hence there is no obvious
relation between symplectic immersions of discs
in the unit ball and in the domain $\Omega$.

In general, the existence of an immersed symplectic
disc in $C$ whose boundary is a given Reeb orbit $\gamma$ on $\Sigma$
does not seem to be known. However, we observe in  
Section 4 that such an immersed
symplectic disc always exists if
$\Sigma$ is the boundary of a strictly 
convex domain in $\mathbb{R}^4$.

There is a second numerical invariant for 
a periodic Reeb orbit $\gamma$ on $\Sigma$, the so-called 
\emph{Maslov index} $\mu(\gamma)$.
If $\Sigma$ is the boundary of a compact strictly convex
body $C\subset \mathbb{R}^4$ then the Maslov index
of any periodic Reeb orbit on $\Sigma$ 
is at least three \cite{HWZ98}.
The \emph{action} of $\gamma$ is
defined to be $\int_\gamma\lambda>0$, and the
orbit is \emph{minimal} if its action is
minimal among the actions of all periodic orbits of
the Reeb flow. Ekeland \cite{ekeland} showed that
the Maslov index of a minimal Reeb orbit on $\Sigma$ equals
precisely three.

Our second result relates the Maslov index to the self-linking number
for periodic Reeb orbits on the boundary of compact strictly
convex bodies with geometric control.

\begin{theorem}\label{thm2}
Let $C$ be a compact strictly 
convex body with smooth boundary $\Sigma$. 
If the principal curvatures $a\geq b\geq c$ of $\Sigma$ satisfy the
pointwise pinching condition $a\leq b+c$ then a periodic
Reeb orbit $\gamma$ on $\Sigma$ of Maslov index $3$ 
bounds an embedded symplectic disc in $C$. In particular, 
the self-linking number of $\gamma$ 
equals $-1$. 
\end{theorem}

As an immediate consequence, if $\Sigma$ is as in Theorem \ref{thm2}
then a periodic Reeb orbit $\gamma$ 
on $\Sigma$ of Maslov index 3 is a slice knot
in $\Sigma$. In fact, 
with some additional effort it is possible to show that 
such an orbit is unknotted \cite{H07}.

The proofs of these results use mainly tools from
differential topology and differential geometry. 
In Section 2 we begin to investigate topological 
properties of transverse knots
on the three-sphere $S^3$.
We define a self-intersection
number for a (not necessarily immersed) disc in the closed unit
ball $B^4$ with boundary $\gamma$ which does not have self-intersections 
near the boundary and relate this 
self-intersection number to 
the self-linking number of $\gamma$.
In Section 3 we study topological invariants of immersed discs
in $B^4$ with boundary $\gamma$ and show
Theorem \ref{thm1} and the corollary. 
In Section 4 we look at boundaries of 
strictly convex bodies in $\mathbb{R}^4$ and derive 
Theorem \ref{thm2}.

\section{Self-intersection of surfaces}

In this section we investigate topological invariants of smooth
maps from an oriented bordered surface $S$ 
with connected boundary $\partial S$ 
into an arbitrary
smooth oriented simply connected $4$-dimensional manifold $W$
(without boundary) 
whose restrictions to a
neighborhood of $\partial S$ are embeddings. For the main
application, $W=\mathbb{C}^2=\mathbb{R}^4$.
We use this discussion to investigate 
maps from the closed
unit disc $D\subset \mathbb{C}$ into the compact
unit ball $B^4\subset \mathbb{C}^2$.
For maps which map the oriented boundary $\partial D$ of $D$
to a canonically oriented transverse knot $\gamma$ on $\Sigma$ 
we define a self-intersection number and relate this to 
the self-linking number
of $\gamma$.

Let for the moment $S$ be any compact oriented surface
with connected boundary $\partial S=S^1$.

\begin{definition}\label{boundary regular}
A smooth map
$f:S\to W$, i.e. a map
which is smooth up to and including the
boundary, is called
\emph{boundary regular} if the singular points of
$f$ are contained in the interior of $S$, i.e. if there is a
neighborhood $A$ of $\partial S$ in $S$ such that the
restriction of $f$ to $f^{-1}(f(A))$ is an embedding.
\end{definition}

McDuff investigated in \cite{MD91}
boundary regular \emph{pseudo-holomorphic} discs in
\emph{almost complex $4$-manifolds} $(W,J)$. By definition,
such a pseudo-holomorphic disc is a
smooth boundary regular
map $f$ from the closed
unit disc $D\subset \mathbb{C}$ into $W$ whose differential
is complex linear with respect to the complex structure
on $D$
and the almost complex structure $J$. She defined a
topological invariant for such boundary regular
pseudo-holomorphic discs which depends on
a trivialization of the normal bundle over
the boundary circle.

Our first goal is to find a purely 
topological analog of this construction. For this we say that
two boundary regular maps $f,g:S\to W$
\emph{are contained in the same boundary class} if $g$ coincides with
$f$ near the boundary and is homotopic to $f$ with fixed
boundary. This means that there is a homotopy $h:[0,1]\times
S\to W$ connecting $h_0=f$ to $h_1=g$ with $h(s,z)=f(z)$ for
all $s\in [0,1]$, all $z\in \partial S$. We do not require that
each of the maps $h_s:z\to h_s(z)=h(s,z)$,
$s\in [0,1]$, is boundary regular.
In particular, if $\pi_2(W)=0$ then any two boundary regular
maps $f,g:S\to W$ which coincide near the
boundary $\partial S$ of $S$ are contained in
the same boundary class (recall that we require that
$W$ is simply connected).

There is also the following stronger notion of homotopy
for boundary regular maps.

\begin{definition}\label{homotopy}
A homotopy $h:[0,1]\times S\to W$
is called \emph{boundary regular} if for each $s$
the map $h_s$ is boundary regular and
coincides with $h_0$ near $\partial S$.
\end{definition}

The set of boundary regular
maps in the boundary class of a
map $f:S\to W$ can naturally be
partitioned into boundary regular homotopy classes.

A boundary regular map $f:S\to W$
is an embedding near $\partial S$.
Since $S$ is oriented by assumption, the normal bundle
$L$ of $f(S)$ over the embedded circle $f(\partial S)$ is an 
oriented real two-dimensional subbundle of 
$TW\vert f(\partial S)$.

For each trivialization
$\rho$ of this normal bundle, the self-intersection
number ${\rm Int}(f,\rho)\in \mathbb{Z}$ is defined
as follows \cite{MD91}.
Let $\overline{N}$ be a closed tubular neighborhood of
$f(\partial S)=\gamma$ in $W$ 
with smooth boundary $\partial N$ 
such that $f(S)\cap \overline{N}$ is
an embedded closed annulus $A$ which intersects
$\partial N$ transversely. Let $E\subset W$ be
an embedded submanifold with boundary which contains $A$ and is
diffeomorphic to an open disc bundle over $A$. One of the 
two connected components $(\partial E)_0$ of the 
boundary $\partial E$ of $E$
has a natural identification
with the total space of the normal bundle
$L$ of $f(S)$ over $\gamma$. Remove
$\overline{N}-E$ from $W$ and glue to the boundary $(\partial E)_0$
of the resulting manifold the oriented real two-dimensional
vector bundle $D\times \mathbb{C}\to D$ in such a way that
$\partial D\times \{0\}$ is identified with the curve 
$f(\partial S)=\gamma\subset (\partial E)_0$
% We glue an oriented disc
%$D^\prime$ to $f(S)$
%along the boundary in such a way that the tangent bundle of this
%disc together with the tangent bundle of $f(S)$ near
%$f(\partial S)$
%combines to the tangent bundle of the closed oriented surface
%$S_0$ which we obtain from $S$ by closing the boundary in the
%usual way. We then
%glue the trivial complex line bundle $L^\prime$
%over $D^\prime$ to
%the bundle $L$ along $\rho$, i.e. we
%identify a section of $L^\prime$ over
%$\partial D^\prime$ which extends to
%a trivialization of $L^\prime$ over $D^\prime$
%with the section of $L$ over $f(\partial S)$
%defined by the trivialization $\rho$.
%The \emph{self-intersection number} ${\rm Int}(f,\rho)$ is then
%defined by ${\rm Int}(f,\rho)=c(\rho)(S_0)-2$ where
%$c(\rho)$ is the
%first Chern class of the resulting $2$-dimensional complex
%vector bundle over $S_0$.
and that the fibres $\{x\}\times\mathbb{C}$
$(x\in \partial D)$ match up
with the trivialized normal bundle $L\vert f(\partial S)$ of
$f(S)$ over $\gamma$.
Up to diffeomorphism,
the resulting $4$-dimensional smooth manifold
$W_\rho$ only depends on the homotopy class of the
trivialization $\rho$ and of the boundary class of $f$.
Let $S_0$ be the closed oriented surface obtained by glueing a disc
to the boundary of $S$ in the usual way. The map $f$ naturally
extends to a map $f_0$ of $S_0$ into $W_\rho$.
The \emph{self-intersection number} ${\rm Int}(f,\rho)$ is then
defined to be
the topological self-intersection number
of $f_0(S_0)$ in $W_\rho$. Thus
${\rm Int}(f,\rho)$
is the number of intersections of $f(S)$ with a surface $f^\prime$
which is a generic perturbation of $f(S)$ and such that $f(\partial S)$
is pushed into the direction given by $\rho$.

In the next lemma we determine the boundary regular homotopy
classes in a fixed boundary class.

\bigskip

\begin{lemma}\label{self}
Let $f:S\to W$ be a smooth boundary regular map.
Choose a trivialization
$\rho$ of the oriented normal
bundle of $f(S)$ over $f(\partial S)$. Then the
assignment which associates to a boundary regular homotopy class
of maps in the boundary class of $f$ its self-intersection number
with respect to $\rho$ is a bijection onto ${\rm Int}(f,\rho)+2\mathbb{Z}$.
Moreover, if $f$ is an embedding then each such class
can be represented by an embedding.
\end{lemma}
\begin{proof}
Let $f:S\to W$ be
a smooth boundary regular map.
Write $\gamma=f(\partial S)$ and let $u:S\to W$ be a
smooth boundary regular
map in the boundary class of $f$. This means that there is
a homotopy $h:[0,1]\times S\to W$
connecting $h_0=f$ to $h_1=u$ with fixed boundary.
The maps $u,f$ are contained in the same
boundary regular homotopy class if and only if
this homotopy can be chosen in such a way that there is a tubular
neighborhood $N$ of $\gamma$ such that the intersection of
$h_s(S)$ with $N$ is independent of $s$.

Choose such an open tubular neighborhood $N$ of $\gamma$ with smooth
boundary $\partial N$ which is sufficiently small that both $f(S)$
and $u(S)$ intersect $N$ in a smooth annulus containing $\gamma$ as
one of its two boundary components. We may assume that there is a
compact subsurface $C\subset S$ with smooth boundary 
$\partial C$ such that $S-C$
is an annulus neighborhood of $\partial S$ and that
$f(S-C)=u(S-C)=f(S)\cap N=u(S)\cap N$. Then $u\vert _C$ and $f\vert
_C$ can be combined to a map into $W-N$ of the closed oriented
surface $\tilde S$ which we obtain from $C$ by gluing two copies of
$C$ along the boundary with an orientation reversing boundary
identification. This map is homotopic in
$W-N$ to a constant map if and only if $u$ and $f$ are contained in the
same boundary regular homotopy class.

Now $W$ is simply connected by assumption and
$N$ is homeomorphic to a $3$-ball-bundle over a circle,
with boundary $\partial N\sim \gamma\times S^2$. Thus by
van Kampen's theorem, $W-N$ is simply connected and
the second homotopy group $\pi_2(W-N)$ coincides with the
second homology group $H_2(W-N,\mathbb{Z})$ via
the Hurewicz isomorphism. Since two boundary regular maps
in the same boundary class are homotopic with fixed boundary, we
conclude that the family of
boundary regular homotopy classes of maps in the boundary class of
$f$ can be identified with the kernel of the
natural homomorphism
$H_2(W-N,\mathbb{Z})\to H_2(W,\mathbb{Z})$.

To compute this group,
we use the long exact homology sequence of the pair $(W,W-N)$
given by
\begin{align}
\dots  \rightarrow H_3 & (W,\mathbb{Z}) &\rightarrow
H_3 & (W,W-N,\mathbb{Z}) &\rightarrow H_2 &(W-N,\mathbb{Z}) & \rightarrow
H_2 &(W,\mathbb{Z}) & \rightarrow \dots \notag
\end{align}
Excision shows that $H_3(W,W-N,\mathbb{Z})=
H_3(\overline{N},\partial N,\mathbb{Z})$ where $\overline{N}$ is the
closure of $N$. Since $\overline{N}=\gamma\times B^3=
S^1\times B^3$ where $B^3$
denotes the closed unit ball in $\mathbb{R}^3$, the group
$H_3(W,W-N,\mathbb{Z})$ is cyclic and generated by a ball
$\{z\}\times (B^3,S^2)$ where $z\in \gamma$ is any fixed point. 

Every
singular homology class $v\in H_3(W,\mathbb{Z})$ can be represented
by a piecewise smooth singular cycle $\sigma$ whose image
is nowhere dense in $W$. On the other hand, the curve
$\gamma$ is contractible in $W$ and therefore there is a smooth
isotopy of $W$ which moves $\sigma$ away from $N$. Thus 
the image of
$H_3(W,\mathbb{Z})$ under the natural homomorphism $H_3(W,\mathbb{Z})\to
H_3(W,W-N,\mathbb{Z})= H_3(\overline{N},\partial N,\mathbb{Z})$
vanishes and hence by exactness, the kernel of the natural
homomorphism
$H_2(W-N,\mathbb{Z})\to H_2(W,\mathbb{Z})$ is isomorphic to
$\mathbb{Z}$ and generated by a sphere $e=\{z\}\times S^2\sim 1\in
\pi_2(\partial N)=\mathbb{Z}$ for some $z\in \gamma$. As a
consequence, every boundary regular homotopy class in the boundary
class of $f$ can uniquely be represented in the form $[f]+ke$ where
$[f]$ denotes the boundary regular homotopy class of $f$ and where
$k\in \mathbb{Z}$.

Next we show that if $f$ is an embedding then each of these classes
can be represented by an embedding as well. For this note that after
possibly replacing $N$ by a smaller tubular neighborhood of $\gamma$
we may assume that
$f^{-1}(\overline{N})$ is a closed annulus neighborhood of $\partial S$
and that  for some $z\in \gamma$, 
the sphere $M=\{z\}\times S^2\subset\partial N$ intersects $f(S)$
transversely in a single point $x$. The orientation of the surface
$S$ defines uniquely an orientation of $M$ such that
$T_xW=T_x(f(S))\oplus T_xM$ as oriented vector spaces. Using
standard surgery near the transverse intersection point $x$ we can
attach the sphere $M$ to the surface $f(S)$ as follows (see
p.38 of \cite{GS99}). There is a closed neighborhood $V$ of $x$ in
$W-\gamma$ which is diffeomorphic to a closed ball and such that the
intersections $f(S)\cap V$, $M\cap V$ are smooth discs which
intersect transversely in the single point $x$. The boundaries of
these discs are two disjoint oriented circles in the boundary
$\partial V\sim S^3$ of $V$. These circles define the Hopf link in
$S^3$ and therefore they form the oriented boundary of a smooth
embedded annulus in $\partial V$. The surgery replaces $(f(S)\cup
M)\cap V$ by such an annulus (which can be
done smoothly). We obtain in this way a compact
oriented bordered surface which can be represented by a boundary
regular map $g:S\to W$ which coincides with $f(S)$ near the
boundary. The surgery does not change relative homology classes
(p.38 of \cite{GS99}) and hence $g(S)$ is homologous to $[f]+e$
via an identification of $M$ with 
a generator $e$ of the kernel
of the natural map $H_2(W-N,\mathbb{Z})\to H_2(W,\mathbb{Z})$. 
In other words, the
embedded surface which we just constructed represents the boundary
regular homotopy class $[f]+e$ in the boundary class of $f$. In the
same way we can also construct a surface which represents the boundary
regular homotopy class $[f]-e$ by attaching to $f$ a sphere equipped
with the reverse orientation. Namely, we also can connect the
boundaries of the discs $f(S)\cap V,M\cap V$ with an embedded
cylinder whose oriented boundary is the union of the oriented
boundary of $f(S)\cap V$ with the boundary of $M\cap V$ equipped
with the reversed orientation. Repeating this procedure finitely
many times with different basepoints we obtain an embedding in every
boundary regular homotopy class of maps in the boundary class of
$f$.

Let $\rho$ be a trivialization
of an oriented normal bundle of $f(S)$ along
$\gamma=f(\partial S)$. We are left with showing that a boundary
regular homotopy class in the boundary class of $f$ is determined
by its self-intersection number with respect to $\rho$. For this
let again
$M=\{z\}\times S^2\subset \partial N$ be an oriented
embedded sphere as above which intersects $f(S)$
transversely in a single point $x$. Assume that the index of
intersection between $f(S)$ and $M$ with respect to the given
orientations is positive.
Let $g:S\to W$ be the map constructed above
with $[g]=[f]+e$.
Using the above notations, it is enough
to show that
${\rm Int}(g,\rho)= {\rm Int}(f,\rho)+2$.
However, this can be seen as follows.

As above, denote by $W_\rho$ the manifold used for the definition of
the self-inter\-sec\-tion number ${\rm Int}(f,\rho)$. Recall that up
to diffeomorphism, the manifold $W_\rho$ only depends on $\rho$ and
the boundary class of $f$. In particular, we may assume that
$W_\rho$ contains the images $\Gamma_f,\Gamma_g$ of the closed
surface $S_0$ under the natural extensions of the maps $f,g$. The
self-intersection numbers of the surfaces $\Gamma_f,\Gamma_g$ in the
manifold $W_\rho$ can now be compared via
\begin{eqnarray*}
{\rm Int}(g,\rho)&=&\Gamma_g\cdot\Gamma_g
=(\Gamma_f+e)\cdot(\Gamma_f+e)\\
&=&\Gamma_f\cdot\Gamma_f+2\Gamma_f\cdot e+e\cdot e
={\rm Int}(f,\rho)+2
\end{eqnarray*}
since the topological self-intersection of the sphere
$e$ in $W_\rho$ vanishes.
But this just means that the assignment which
associates to a boundary regular homotopy class in the
boundary class of $f$ its self-intersection number
with respect to $\rho$ is a bijection onto
${\rm Int}(f,\rho)+2\mathbb{Z}$.
%Assume that $f$ is an immersion and let $L$ be the normal
%bundle of $f(S)$ in $W$. Note that $L$ is an oriented
%real two-plane bundle. Since $S$ is a bordered surface,
%such a bundle admits a global trivialization. Choose
%a trivialization $\zeta$ of $L$ over $\gamma$
%which extends to a global trivialization of
%a complex line-subbundle $L$ of $TW$ over $f(S)$.
%Then ${\rm Int}(f,\rho)$ is just the winding number
%of $\rho$ with respect to $\zeta$. Now let
%$\zeta^\prime$ be the trivialization over
%$\gamma$ of the normal bundle $L^\prime$
%for the surface $g:S\to W$
%which we constructed above from $f$ and $M$ and
%which extends to a global trivialization of $L^\prime$.
%Then the winding number of $\zeta^\prime$ with respect to
%$\zeta$ equals $-1$ and therefore
%${\rm Int}(g,\rho)={\rm Int}(f,\rho)-1$ as claimed.
From this the lemma follows.
\end{proof}

From now on we assume that the $4$-dimensional
manifold $W$ is equipped with a smooth almost
complex structure $J$. 

\begin{definition}\label{infinitesholo}
A smooth boundary regular
map $f:S\to W$ is called \emph{boundary holomorphic}
if for each $z\in \partial S$ the
tangent plane of $f(S)$ at $z$ is a complex line
in $(TW,J)$ whose orientation coincides with
the orientation induced from the orientation of $S$.
\end{definition}

If $f:S\to W$ is boundary regular and boundary holomorphic
then
the pull-back $f^*TW$ under $f$ of the tangent bundle
of $W$ is a 2-dimensional complex vector bundle
over $S$. Since $f$ is
boundary holomorphic, the restriction to $\partial S$ of 
the tangent bundle $TS$ of $S$ 
is naturally a complex line-subbundle
of $f^*TW\vert_{\partial S}$.
Then the normal bundle of $f(S)$ over $\gamma$ 
can be identified with a complex
line subbundle of $(TW,J)\vert f(\partial S)$ as well.
Every trivialization $\rho$ of this normal bundle
defines as before a smooth manifold $W_\rho$. This manifold
admits a natural almost complex structure
extending the almost complex structure on the
complement of a small tubular neighborhood of $f(\partial S)$
in $W$. In particular, if we denote as before by
$f_0$ the natural extension of $f$ to the
closed surface $S_0$ then the pull-back
bundle $f_0^*TW_\rho$ is a complex two-dimensional
vector bundle over $S_0$. 

Up to homotopy, this
bundle only depends on $\rho$ and the boundary
class of $f$. Namely, any homotopy $h_s$ of $f=h_0$ 
which is the identitity near the boundary induces a homotopy of the
pull-back bundles $h_s^*TW_\rho$. Now if $f_0,f_1$ are the 
extensions of $h_0,h_1$ to the closed surface $S_0$ then 
since the homotopy $h_s$ 
is the identity near the boundary, it determines a 
homotopy of the
complex pull-back bundle $f_0^*TW_\rho$ to the complex
pull-back bundle $f_1^*TW_\rho$.
Let $c(\rho)$ be the evaluation on $S_0$ of the first Chern class
of this bundle. 

Changing the trivialization $\rho$ by a full
positive (negative) twist in the group $U(1)\subset
GL(1,\mathbb{C})$ changes both the self-intersection number ${\rm
Int}(f,\rho)$ and the Chern number $c(\rho)$ by $ 1$ $(-1)$
(see \cite{MD91}).
In particular, there is
up to homotopy a unique trivialization $\rho$
of the complex normal bundle
of $f(S)$ over $f(\partial S)$ such that $c(\rho)=2$.
We call such a trivialization a
\emph{preferred} trivialization.
By the above observation,
a preferred trivialization only depends on
the boundary class of $f$
but not on the boundary regular homotopy class of $f$.
Moreover, the complex normal bundle of a boundary
holomorphic boundary regular map $f:S\to W$ only depends on
the oriented boundary circle $f(\partial S)$.  

\begin{definition}\label{selfintersectionnumber}
The \emph{self-intersection number}
${\rm Int}(f)$ of a 
boundary holomorphic boundary regular map $f:S\to W$ 
is the self-intersection number ${\rm Int}(f,\rho)$ of $f$ 
with respect to a preferred
trivialization $\rho$ of the complex
normal bundle of $f(S)$ over $f(\partial S)$.
\end{definition}

Lemma \ref{self} implies

\begin{kor}\label{uniquedetermined}
Let $f:S\to W$ be
boundary regular and boundary holomorphic.
Then a boundary regular homotopy
class in the boundary class of $f$ is uniquely determined
by its self-intersection number.
\end{kor}

Now consider the standard unit sphere 
$S^3$ in $\mathbb{C}^2$ which bounds the standard open
unit ball $B^4_0\subset \mathbb{C}^2$.
The \emph{contact distribution} on 
$S^3$ is the unique smooth two-dimensional
subbundle $\xi$ of $TS^3$ which is invariant under the
(integrable) complex structure $J$ on $\mathbb{C}^2$. 
A smooth embedding $\gamma:S^1\to S^3$ is 
\emph{transverse} if its tangent 
$\gamma^\prime$ is everywhere transverse to $\xi$.
Let $N$ be the outer normal field of $S^3$. Then
$JN$ is tangent to $S^3$ and orthogonal to $\xi$.
The tangent of the transverse knot $\gamma$ can be written
in the form 
\[\gamma^\prime(t)=a(t)JN(\gamma(t))+B(t)\]
where $B(t)\in \xi$ for all $t$ and where $a(t)\not=0$.
Assume that $\gamma$ is oriented in such a way that
$a(t)$ is positive for all $t$. 
If as  in the introduction we denote by $\lambda$ 
be the restriction to $S^3$ of 
the radial one-form $\lambda_0$ on 
$\mathbb{C}^2$ then this 
orientation of $\gamma$ is determined by the requirement that
the evaluation of $\lambda$ on $\gamma^\prime$
is positive, and we call it \emph{canonical}.
If the transverse knot $\gamma$ is canonically oriented
then $J\gamma^\prime$ points
inside the ball $B^4_0$. Thus if $B^4=B^4_0\cup S^3$ 
denotes the closed unit ball
then for every canonically oriented transverse
knot $\gamma$ on $S^3$ there is a
smooth boundary regular boundary holomorphic
map $f:S\to B^4$ with $f(\partial S)=\gamma$ and
$f^{-1}(S^3)=\partial S$ whose restriction to 
a neighborhood of the boundary is an embedding: Just choose 
a smooth embedding of a closed annulus $A$ into
$B^4$ which maps one of the boundary circles $\zeta$ of $A$ diffeomorphically
onto $\gamma$ and whose tangent plane at a point in $\zeta$ is 
$J$-invariant. In particular, $A$ meets $S^3$
transversely along $\gamma$ and hence we may assume that
$A\cap S^3=\gamma$.  
Extend this embedding in an arbitrary way 
to a smooth map of the surface $S$ (with 
the annulus $A$ as a neighborhood of $\partial S$) into $B^4$
which is
always possible since $B^4$ is contractible.

Define a boundary regular map $f:(S,\partial S)\to (B^4,S^3)$
to be \emph{boundary transverse} if $f$ is transverse to 
$S^3$ along the boundary.
Since $B^4$ is contractible, the above observation implies that
every boundary regular boundary transverse map 
$f:(S,\partial S)\to (B^4,\gamma)$
can be homotoped within the
family of such maps to a boundary regular boundary holomorphic
map $f^\prime:(S,\partial S)\to (B^4,\gamma)$. 
In particular, the oriented normal bundle of $f$ over $\gamma$ is
naturally homotopic to the oriented normal bundle of
$f^\prime$ over $\gamma$. The map $f^\prime$ is used to
calculate the preferred trivialization of this normal bundle.
Then the self-intersection number ${\rm Int}(f)$ can be defined
as the self-intersection number of $f$ with respect to the induced
trivialization of the oriented normal bundle of $f$. This self-intersection
number coincides with the self-intersection number ${\rm Int}(f^\prime)$ 
of $f^\prime$ and by Corollary \ref{uniquedetermined}, it
only depends on $f$.

More generally, let $\Sigma$ be the boundary of a bounded
domain $\Omega\subset \mathbb{C}^2$ 
which contains the origin $0$ in its interior and
which is star-shaped with respect
to $0$. The contact form is the 
restriction $\lambda$ to $\Sigma$ of the
radial one-form $\lambda_0$ on
$\mathbb{C}^2$.

Let $N$ be the outer normal field of $\Sigma\subset
\mathbb{C}^2$. Since the one-form $\lambda_0$ can
also be written in the form
$(\lambda_0)_p(Y)=\frac{1}{2}\langle Jp,Y\rangle$
$(p\in \mathbb{C}^2,Y\in T_p\mathbb{C}^2$ and where $\langle,\rangle$ is
the euclidean inner product),
the \emph{Reeb vector field} $X$ on
$\Sigma$ is given by
\[X(p)=\phi(p)JN(p)\]
where
\[\phi(p)=\frac{2}{\langle p,N(p)\rangle}>0.\]
Namely, for $p\in \Sigma$ we have
\[d\lambda_p(X,\cdot)=\phi(p)\omega_0(JN(p),\cdot)=
-\phi(p)\langle N(p),\cdot\rangle=0\]
on $T_p\Sigma$ and
\[\lambda_p(X)=\frac{1}{2}\langle Jp,X\rangle=
\frac{1}{2}\phi(p)\langle Jp,JN(p)\rangle =1.\]
In particular, if $\gamma$ is a Reeb orbit on $\Sigma$
then $J\gamma^\prime$ is transverse to $\Sigma$ and points inside
the domain
$\Omega$. As a consequence, as for transverse knots
on $S^3$, if we denote
by $C=\Omega\cup \Sigma$ the closure
of $\Omega$ then for every boundary regular boundary transverse map 
$f:S\to C$ whose boundary $f(\partial S)$
is a periodic Reeb orbit on $\Sigma$, the self-intersection
number ${\rm Int}(f)$ of $f$ is well defined.

The next lemma shows that in both cases, the self-intersection
number of such a map $f:(S,\partial S)\to (C,\gamma)$ only
depends on $\gamma$. For a convenient formulation,
let $C$ be the closure
of a bounded star-shaped domain in $\mathbb{C}^2$ and
call a smoothly embedded closed curve $\gamma$ in the
boundary $\Sigma$ of $C$   
\emph{admissible} if either $C$ is the unit ball, $\Sigma=S^3$ and
$\gamma$ is a canonically oriented transverse knot or
if $\gamma$ is a Reeb orbit on $\Sigma$.

\begin{lemma}\label{homotopyinv}
Let $\gamma\subset \Sigma$ be an admissible curve. Then
any two boundary regular boundary transverse
maps $f:(S,\partial S)
\to (C,\gamma),g:(S^\prime,\partial S^\prime)\to (C,\gamma)$ 
have the same self-intersection number.
\end{lemma}
\begin{proof}
Let $f:S\to C,g:S^\prime\to C$ be any two 
boundary regular boundary transverse maps with boundary 
an admissible closed curve $\gamma$. 
The inner normals
of the surfaces $f(S),g(S^\prime)$
along $\gamma=f(\partial S)=g(\partial S^\prime)$
point strictly inside the domain $C$.
After a small deformation through
boundary regular boundary transverse maps we
may assume that there is a small annular
neighborhood $A$ of the boundary of $S$,
an annular neighborhood $A^\prime$ of the
boundary of $S^\prime$ and a homeomorphism
$\phi:A\to A^\prime$ which maps $\partial S$
to $\partial S^\prime$ and is such that
$g(\phi(x))=f(x)$ for all $x\in A$. We may moreover
assume that
the restrictions of
$f,g$ to $f^{-1}(A),g^{-1}(A^\prime)$ are embeddings.
Since $f,g$ are boundary regular, after possibly
modifying $f,g$ once more with a small
boundary regular homotopy which pushes interior intersection
points of $f(S),g(S^\prime)$ with $\Sigma$ into
the interior $\Omega$ of $C$ we may assume that
there is a compact star-shaped set
$K\subset \Omega$ such that
$f(S-A)\subset K,g(S^\prime-A)\subset K$. But then
the restrictions of $f,g$ to $S-A,S^\prime-A^\prime$
are maps of surfaces $S-A,S^\prime-A^\prime$
into $K$ with the same boundary curve $\gamma^\prime$.
Now $K$ is contractible and hence the maps $f\vert S-A,
g\vert S^\prime-A^\prime$ 
define the same relative
homology class in $H_2(K,\gamma^\prime;\mathbb{Z})$.
By Lemma 2.3 and its proof, this implies that 
the self-intersection numbers of
$f,g$ indeed coincide.
\end{proof}

As a consequence, we can define:

\begin{definition}\label{self-intersection}
Let $\gamma\subset \Sigma$ be 
an admissible curve.
The \emph{self-intersection number}
${\rm Int}(\gamma)$ of $\gamma$
is the self-intersection number of a boundary regular
boundary transverse 
map $f:S\to C$ with
boundary $f(\partial S)=\gamma$.
\end{definition}

The final goal of this section is to calculate the self-intersection 
number of an admissible curve $\gamma$ on
$\Sigma$.
For this we begin with calculating the preferred trivialization of the
normal bundle of the complex line subbundle of $T\mathbb{C}^2\vert \gamma$
spanned by the tangent $\gamma^\prime$ of $\gamma$. 
Note that this normal bundle can naturally be
identified with the restriction of the contact bundle $\xi$ to $\gamma$.

To this end let $D$ be the closed unit disc in $\mathbb{C}$ and let 
$f:D\to C$ be a boundary holomorphic boundary regular
\emph{immersion} which maps $\partial D$ diffeomorphically
onto the canonically oriented admissible curve $\gamma$. 
The image under $df$ of the inner normal
of $D$ along $\partial D$
points strictly inside $C$. Let
\[\hat M:(z_1,z_2)\to (-\overline{z}_2,\overline{z}_1)\]
be a $J$-orthogonal $\langle\cdot ,\cdot\rangle$-compatible
almost complex structure on $\mathbb{C}^2$
where as usual, $z\to \overline{z}$ is complex conjugation.
Then we obtain a
trivialization of the complex
vector bundle $(f^*T\mathbb{C}^2,J)$ over $D$ by the
sections $X_1=df(\frac{\partial}{\partial x}),
X_2=\hat Mdf(\frac{\partial}{\partial x})$ (with a slight
abuse of notation).

The trivialization of the tangent bundle $df(TD)\vert \gamma$
of $f(D)$ over $\gamma$ defined by the tangent
$\gamma^\prime$ of the admissible curve $\gamma$
has rotation number one with respect to the
trivialization $df(\frac{\partial}{\partial x})$.
Since the tangent bundle of the two-sphere
$S^2$ has Chern number $2$
and is obtained by glueing the tangent bundles
of two standard discs $D_1,D_2$ along the 
boundary using the trivializations defined by the 
tangent field of the boundary, the restriction of 
the section $X_2$ to $\gamma$ defines the preferred
trivialization of the complex normal bundle $L$ 
over $\gamma$. Now $\hat M$ is complex anti-linear
and therefore the trivialization of $L$ over $\gamma$
defined by the section $\hat M\circ \gamma^\prime$ has
rotation number $-1$ with respect to the preferred 
trivialization. Recall that the
preferred trivialization of the complex normal bundle
of $\gamma$ only depends on the boundary
class of an infinitesimally holomorphic immersion of a
surface $S$ into $\mathbb{C}^2$.

To the admissible curve $\gamma$
we can also associate its \emph{self-linking number}
${\rm lk}(\gamma)$ (see \cite{eliash} and the introduction).
The following proposition relates these two numbers.

\begin{proposition}\label{inter-link}
Let $\gamma$ be an admissible curve on
$\Sigma$.
Then the self-intersection number ${\rm Int}(\gamma)$ of
$\gamma$ equals ${\rm lk}(\gamma)+1$.
\end{proposition}
\begin{proof}
\begin{figure}[ht]
\psfrag{N}{$N(p)$}
\psfrag{F}{$F$}
\psfrag{F2}{$F^\prime$}
\psfrag{S}{$\Sigma$}
\psfrag{N2}{$\nu$}
\psfrag{N1}{$\nu^\prime$}
  \centering
  \includegraphics[scale=0.35]{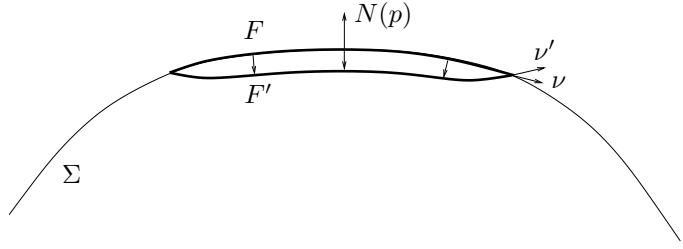}
  \caption{The surfaces $F$ and $F^\prime$}
  \label{seif}
\end{figure}
Let $N$ be the outer normal field of $\Sigma$ and let
$L$ be the \emph{complex subbundle} of $T\Sigma$,
i.e. the $2$-dimensional subbundle which is 
invariant under the complex structure $J$. 
If $\Sigma=S^3$ then this is just the contact bundle.
The image of
the outer normal $N$ of $\Sigma$ 
under the $J$-orthogonal
$\langle,\rangle$-compatible almost complex structure
$\hat M$ is a global section of the bundle $L$.
Let $F\subset \Sigma$ be a \emph{Seifert surface}
for the admissible curve $\gamma$, i.e. $F$
is an embedded oriented bordered surface in $\Sigma$
with boundary $\gamma$.
Let $N_F$ be the oriented normal field
of $F$ in $\Sigma$ with respect to the restriction of the
euclidean metric $\langle,\rangle$.
Since $\gamma$ is transverse to $L$ we may assume that
for every $x\in \gamma$ the vector $N_F(x)$ is contained in
the fibre $L_{x}$ at $x$ of
the complex line bundle $L$.
The self-linking number of $\gamma$ is
therefore the winding number of the section
$x\to M(x)=\hat MN(x)$ of $L\vert \gamma$
with respect to the trivialization of $L\vert \gamma$
defined by the section $x\to N_F(x)$.

Let $F^\prime\subset C$ 
be the embedded surface which we obtain by pushing $F$ 
slightly in the direction $-N$ as in
Figure \ref{seif}. Then $F^\prime$ is an
embedded surface in $C$ which is 
boundary regular and boundary transverse. The
restriction of $N_F$ to $\gamma$ extends to a global trivialization
of the oriented normal bundle of the surface $F$ and hence
$F^\prime$ in $\mathbb{C}^2-\gamma$. Thus the self-intersection
number of $F^\prime$ with respect to the trivialization defined by
$N_F$ vanishes. Since by the above observation the winding number
of the section $M$ of $L\vert \gamma$ with respect to the preferred
trivialization of $L\vert \gamma$ equals $-1$, the self-intersection
number ${\rm Int}(\gamma)$ equals the winding number of $M$ with
respect to the trivialization of $L\vert \gamma$ 
defined by $N_F$ plus one. This shows
the proposition.
\end{proof}

\section{Topological invariants of immersed discs}

As in Section 2, we denote by $S$ a compact oriented
surface with connected boundary $\partial S=S^1$. Let $(W,J)$
be a smooth simply connected 4-dimensional manifold equipped with
a smooth almost complex structure $J$. In this section we
investigate topological invariants of boundary
holomorphic boundary regular \emph{immersions} $f:S\to (W,J)$.
For this we use the assumptions and notations from
Section 2. Recall in particular the definition of the
self-intersection number ${\rm Int}(f)$ of $f$.

Let $G(2,4)$ be the Grassmannian of oriented (real)
2-planes in $\mathbb{R}^4=\mathbb{C}^2$. This Grassmannian is just
the homogeneous space $G(2,4)=SO(4)/SO(2)\times SO(2)=
S^2\times S^2$, in particular the second homotopy group of $G(2,4)$
coincides with its second homology group and is isomorphic to
$\mathbb{Z}\oplus\mathbb{Z}$.
Moreover, there are generators
$\tilde e_1,\tilde e_2$ of $H_2(G(2,4);\mathbb{Z})$
such that with respect to these generators,
the homological intersection form $\iota$ is the
symmetric form represented by the $(2,2)$-matrix $(a_{i,j})$ with
$a_{1,1}=a_{2,2}=0$ and $a_{1,2}=a_{2,1}=\iota(\tilde e_1,\tilde e_2)=1$.

The complex projective line $\mathbb{C}P^1=S^2$ of 
\emph{complex} oriented
lines in $\mathbb{C}^2$ for the standard complex structure $J$ 
is naturally embedded in $G(2,4)$. 
Its homotopy class is
the generator of an infinite cyclic subgroup $Z_1$ of
$\pi_2(G(2,4))$. We call this generator the \emph{canonical
generator} of $Z_1$ and denote it by $e_1$. 
The anti-holomorphic sphere of all complex oriented
lines for the complex structure $-J$ 
is homotopic and hence homologous in $G(2,4)$ to
the complex projective line $\mathbb{C}P^1$.
Namely, the complex structures $J,-J$ define the 
same orientation on $\mathbb{R}^4$ 
and hence $J$ can be connected to $-J$ by a 
continuous
curve of linear complex structures on $\mathbb{R}^4$. 
This curve then determines a homotopy of $\mathbb{C}P^1$ 
onto the anti-holomorphic 
sphere of complex oriented lines for $-J$. Since these
two spheres are disjoint, the self-intersection number of 
the class in $H_2(G(2,4);\mathbb{Z})$ defined by $e_1$ vanishes.

A second
infinite cyclic subgroup $Z_2$ 
of $\pi_2(G(2,4))$ is defined as follows. Let
$S^2\subset \mathbb{R}^3\subset \mathbb{R}^4$ be the standard unit
sphere. The map which associates to a point $y\in S^2$ the
oriented tangent plane of $S^2$ at $y$, viewed
as a 2-dimensional oriented
linear subspace of $\mathbb{R}^4$, defines a smooth map of
$S^2$ into $G(2,4)$. Its homotopy class 
$e_2$ generates a subgroup $Z_2$ of
$\pi_2(G(2,4))$. We call $e_2$ the canonical
generator of $Z_2$. 
Now the tangent bundle of $S^2\subset \mathbb{R}^3 \subset \mathbb{C}^2$
intersects the complex projective line
$\mathbb{C}P^1\subset G(2,4)$ 
in precisely one point (which is the
tangent space of $S^2$ at $(0,0,1,0)$). This intersection is
transverse with positive intersection index. Therefore we have
$\iota(e_1,e_2)=1$ and hence the elements
$e_1,e_2$ generate $\pi_2(G(2,4))$.

Let $W$ be a simply connected $4$-dimensional
manifold with smooth
almost complex structure $J$. Equip
the tangent bundle $TW$ of $W$ with a
$J$-invariant Riemannian metric $\langle ,\rangle$. Let ${\cal
G}\to W$ be the smooth fibre bundle over $W$ whose fibre at a
point $x\in W$ consists of the Grassmannian of oriented 2-planes
in $T_x W$. Let $S$ be a compact oriented
surface with
connected boundary $\partial S\sim S^1$. Then every smooth
boundary regular \emph{immersion} $f:S\to W$ defines a smooth map
$Gf$ of $S$ into the bundle
${\cal G}$ by assigning to a point $x\in S$ the oriented
tangent space $df(T_xS)$ of $f(S)$ at $f(x)$.
The complex pull-back bundle $(f^*TW,J)$ over $S$ admits
a \emph{complex} trivialization. This
trivialization can be chosen to be of the form $df(X),V$ where $X$
is a global nowhere vanishing section of the tangent bundle
$TS$ of $S$ and $V$ is a global
section of the $\langle, \rangle$-orthogonal complement of the
complex line subbundle of $f^*TW$ which is spanned by the
section $df(X)$. With respect to this complex trivialization of
$f^*TW$, the
pull-back $f^*{\cal G}$ of the bundle ${\cal G}$ can
naturally be represented as a product $S\times G(2,4)$. If $f$ is
boundary holomorphic, i.e. if for every
$z\in \partial S$ the tangent space $df(T_zS)\subset
T_{f(z)}W$ is $J$-invariant
and if moreover its orientation coincides with the orientation
induced by $J$, then in the above
identification of $f^*{\cal G}$ with $S\times G(2,4)$
the circle of tangent planes of $f(S)$ over $f(\partial S)$
is given by the curve $\partial S\times L_0$ in $\partial S\times
G(2,4)$ where $L_0=\mathbb{C}\times \{0\}\subset \mathbb{C}^2$ is
a fixed complex line. Thus in this case the map $Gf$ can be
viewed as a smooth map of the surface $S$ into the Grassmannian
$G(2,4)$ which maps the boundary $\partial S$ of $S$ 
to the single complex line $L_0$.
In other words, if we let $\tilde S$ be the closed
oriented surface obtained from $S$ by collapsing
$\partial S$ to a point then 
$Gf$ defines a map of $\tilde S$ 
into $G(2,4)$. This map
then defines a homotopy class of maps $\tilde S\to G(2,4)$ and a
homology class $[Gf]\in H_2(G(2,4),\mathbb{Z})$.

The following definition strengthens Definition \ref{homotopy}.

\begin{definition}\label{reghomotop}
A smooth homotopy $h:[0,1]\times S\to W$ is \emph{regular}
if $h$ is boundary regular and if moreover for
every $s\in [0,1]$ the map $h_s$ is a boundary
holomorphic immersion.
\end{definition}

If $f,g:S\to W$ are two boundary regular boundary
holomorphic
immersions which are regularly homotopic, i.e. which
can be connected by a regular homotopy,
then the maps $Gf$ and $Gg$ are homotopic.
Namely, if $h:[0,1]\times S\to W$ is a regular homotopy
connecting $h_0=f$ to $h_1=g$, then there is a complex
trivialization of the complex pull-back bundle
$(h^*TW,J)$ over $[0,1]\times S$ whose restriction to
$[0,1]\times
\partial S$ does not depend on $s\in [0,1]$ and is determined
as before by the section $dh_s(X)$ where $X$ is a global 
nowhere vanishing section of $TS$.
This trivialization then defines
an identification of the bundle $h^*{\cal G}$ with
$[0,1]\times S\times G(2,4)$.
For each $s\in [0,1]$ the tangent planes of the immersion $h_s$
define a smooth section of the bundle $h^*{\cal G}$ over
$\{s\}\times S$ and hence a smooth map of $S$ into
$G(2,4)$. This map depends smoothly on $s$ and
maps the boundary $\partial S$ of $S$ to a single point.
Thus by continuity, the homotopy class of the tangent map of $h_s$ is
independent of $s\in [0,1]$ and hence it is an invariant of
regular homotopy.

There is another way to obtain an invariant of regular homotopy.

\begin{definition}\label{tangentialindex}
The \emph{tangential index}
${\rm tan}(f)$ of a boundary regular
immersion $f$ whose only 
self-intersection points are transverse double points
is defined to be the number of  
self-intersection points of $f$ counted with signs.
\end{definition}

If a boundary regular immersion
$f:S\to W$ has self-intersection points which are not
transverse double points then it can be perturbed with a regular
homotopy to an immersion whose only self-intersection
points are transverse double points and
whose tangential index is independent of the
perturbation (see e.g. \cite{MD91}). 
The tangential index is 
invariant under regular homotopy.

If we consider more specifically boundary regular
boundary holomorphic immersions of \emph{discs}
then we can derive a more precise result. For
this recall from Section 2 that for every boundary regular
boundary holomorphic map $f:D\to W$ there
is a preferred trivialization of the normal bundle of
$f(D)$ over $f(\partial D)$.
On the other hand, there is a trivialization
$N$ of the oriented normal bundle of $f(D)$ over $f(\partial D)$ which
extends to a global trivialization of the oriented normal bundle
of $f(D)$ in $TW$.

\begin{definition}
The \emph{winding number} ${\rm wind}(f)$ of a boundary
regular boundary holomorphic immersion
$f:D\to W$ is the winding number of the preferred
trivialization of the normal bundle of $f(D)$ over
$f(\partial D)$ with respect to a trivialization which
extends to a global trivialization of the normal bundle of
$f(D)$.
\end{definition}

For the formulation of the following version of  
the well known \emph{adjunction formula} for immersed
boundary regular boundary holomorphic discs, denote  
for a boundary regular boundary holomorphic immersion 
$f:D\to W$ by
${\cal C}_2(Gf)$ the component of $[Gf]$ in the subgroup 
$Z_2$ of $H_2(G(2,4),\mathbb{Z})$, viewed as an integer.

\begin{proposition}\label{intofdiscs}
For a
boundary regular boundary holomorphic
immersion $f:D\to W$ we have
${\rm Int}(f)={\rm wind}(f)+2{\rm tan}(f)$, and
${\rm wind}(f)=2{\cal C}_2(Gf)$.
\end{proposition}
\begin{proof}
Let $f:D\to W$ be a boundary regular
boundary holomorphic immersion. As in Section 2, let $\rho$
be the preferred trivialization of the complex normal bundle
of $f(D)$ over $f(\partial D)$ and use this trivialization
to extend $f$ to an immersion $f_0$ of the two-sphere $S^2$ into
the almost complex manifold $(W_\rho,\tilde J)$. 
Then ${\rm Int}(f)$ is the self-intersection number 
of $f_0(S^2)$ in $W_\rho$. Since $f_0$ is an
immersion, this self-intersection number just equals
$\chi(N)+2{\rm tan}(f_0)$ where $\chi(N)$ is the Euler number
of the normal bundle of $f_0(S^2)$ in $W_\rho$ and where
${\rm tan}(f_0)={\rm tan}(f)$ is the tangential index
defined above (see e.g. Lemma 4.2 of \cite{MD91}
or simply note that the formula is obvious if $f_0$
is an embedding and follows for immersions with only
transverse double points by
surgery at every double self-intersection point which 
increases the Euler class of the normal bundle by $2$ if
the double point has positive index and decreases it
by $2$ if the double point has negative index). 
By our definition of the winding number ${\rm wind}(f)$ of $f$,
this is just the formula stated in the proposition.

To show that ${\rm wind}(f)=2{\cal C}_2(Gf)$, note first that
we have ${\rm wind}(f)=0$ if $[Gf]\in Z_1$.
Namely, using the above notations,
recall that a preferred trivialization $\rho$
of the normal bundle of the disc $f(D)$ 
over $f(\partial D)=
\gamma$ is determined by the requirement that 
the evaluation of the
first Chern class of the complex tangent bundle 
$(TW_\rho,\tilde J)$ 
of $W_\rho$
on the 2-sphere $f_0(S^2)$ equals two.

The tangent plane map of $f$ can be viewed as a map 
$(D,\partial D)\to G(2,4)$ which maps the boundary
$\partial D$ of $D$ to a single point and hence
factors through a map $F:S^2\to G(2,4)$. 
If $[Gf]\in Z_1$ then since $\pi_2(G(2,4))=H_2(G(2,4),\mathbb{Z})$, 
the map $F$ can be homotoped 
to a map $S^2\to \mathbb{C}P^1$.
By construction, this implies that up to homotopy, the
complex vector bundle $(f_0^*TW,\tilde J)$
decomposes as a direct sum $TS^2\oplus N$ of two complex line bundles.
The first Chern class of $(f_0^*TW,\tilde J)$ is then the sum
of the Chern classes of $TS^2$ and $N$. Therefore by
our normalization, 
the first Chern class of the normal bundle 
$N=f_0^*TW/TS^2\to S^2$
vanishes. As a consequence, 
the bundle $N\to S^2$ is trivial
and hence 
the preferred trivialization of the normal bundle
$N$ over $\gamma$ extends to a global trivialization of 
$N$ over 
$S^2$. This shows that
${\rm wind}(f)=0$ if $[Gf]\in Z_1$.

Arguing as in the proof of Lemma \ref{self},
if $g:D\to W$ is any boundary regular boundary
holomorphic immersion with
${\cal C}_2(Gg)=ke_2$ for some $k\in \mathbb{Z}$ then there
is a boundary regular boundary holomorphic immersion
$f:D\to W$ with $f(\partial D)=g(\partial D)$, 
$[Gf]=[Gg]-ke_2\in Z_1$ and such that 
${\rm Int}(f)={\rm Int}(g)$,  
${\rm tan}(f)={\rm tan}(g)+k$
and ${\rm wind}(f)={\rm wind}(g)-2k$.
Namely, such an immersion $f$ can be constructed as follows.
Choose a point $z\in D$ such that there is a small ball
$V\subset W$ about $g(z)$ which intersects $g(D)$ in an
embedded disc $B$ containing $g(z)$. 
Choose an embedded $2$-sphere $\hat S\subset \partial V$
which intersects $B$ transversely in precisely two points, one
with positive and one with negative intersection index.
The tangent bundle of the sphere is a generator 
$e_2$ of the subgroup $Z_2$ of 
$G(2,4)$.
As in the proof of Lemma \ref{self}, 
attaching the sphere to $g(D)$ with
surgery about the intersection
point with positive intersection index results in a disc $u$ 
which satisfies 
${\rm Int}(u)={\rm Int}(g)$,
$[Gu]=[Gg]+e_2$ and ${\rm tan}(u)=
{\rm tan}(g)-1$. 
Similarly, attaching the sphere to $g(D)$ with surgery
about the intersection point with negative intersection index
results in a disc $u^\prime$ which satisfies 
${\rm Int}(u^\prime)={\rm Int}(g)$,
$[Gu^\prime]=[Gg]-e_2$ and ${\rm tan}(u^\prime)=
{\rm tan}(g)+1$.
From this the proposition is immediate.
\end{proof}

A \emph{complex point} of an immersed
disc $f:D\to W$ is
a point $z\in D$ such that the real two-dimensional subspace
$df(T_zD)$ of $TW$ is invariant under the almost
complex structure $J$. The point is called
\emph{holomorphic} if the orientation of
$df(T_zD)$ induced by the orientation of $D$ coincides
with the orientation induced by the almost complex
structure $J$, and it is called \emph{anti-holomorphic}
otherwise.

\begin{corollary}\label{antiholo}
Let $f:D\to W$ be a boundary regular boundary 
holomorphic immersion. If
$f$ does not have any anti-holomorphic points
then ${\rm Int}(f)=2{\rm tan}(f)$.
\end{corollary}
\begin{proof} Let $f:D\to W$ be a boundary regular boundary 
holomorphic immersion without
any anti-holomorphic point.
Then the tangent map $Gf$ of $f$ does not
intersect the anti-holomorphic sphere of complex lines 
in $\mathbb{C}^2$ equipped with the reverse of the orientation
induced by the complex structure.
Since the anti-holomorphic sphere is homologous in 
$G(2,4)$ to the complex projection line
$\mathbb{C}P^1$ and has vanishing self-intersection
(see the discussion at the beginning of this
section),  
we have $[Gf]\in Z_1$ by consideration of intersection
numbers. The corollary now is immediate from Proposition
\ref{intofdiscs}.
\end{proof}

As in the introduction, denote by $\omega_0$ the standard
symplectic form on $\mathbb{C}^2$. 
An immersion $f:D\to \mathbb{C}^2$ 
is called symplectic if for every $z\in D$ 
the restriction of $f^*\omega_0$
to the tangent plane $T_zD$ does not vanish and defines
the standard orientation of $T_zD$.
As a consequence of Corollary \ref{antiholo} 
we obtain Theorem \ref{thm1} and the corollary
from the introduction.

\begin{corollary}\label{inverse}
Let $\gamma$ either be a transverse knot on the 
standard three-sphere, the boundary of the standard unit ball
$C\subset \mathbb{C}^2$, or a Reeb orbit on the boundary
$\Sigma$ of a domain in $\mathbb{C}^2$ which is 
star-shaped with respect to the origin, with compact closure $C$. 
If $\gamma$ bounds a boundary regular immersed symplectic disc 
$f:(D,\partial D)\to (C,\gamma)$ then ${\rm lk}(\gamma)=
2{\rm tan}(f)-1$.
\end{corollary}
\begin{proof}
By definition, a symplectic immersion does not have
any anti-holomorphic points. Thus 
if $f:(D,\partial D)\to (C,\gamma)$ is a boundary regular boundary
holomorphic immersed symplectic disc then ${\rm lk}(\gamma)=2{\rm tan}(f)-1$
by Proposition \ref{inter-link} and Corollary \ref{antiholo}.

Now if $f:(D,\partial D)\to (C,\gamma)$ is an arbitrary boundary regular
immersed symplectic disc then 
$f$ can be slightly modified with a smooth homotopy to a 
boundary transverse symplectic disc 
without changing the
tangential index since 
being symplectic is an open condition.
Locally near the boundary, this disc can be represented as a graph over
an embedded 
symplectic annulus $A\subset C$ with $\gamma$ as
one of its boundary components whose tangent plane is $J$-invariant
at every point in $\gamma$. 
Now for a fixed nonzero vector
$X\in T\mathbb{C}^2$, the set of all nonzero vectors
$Y\in T\mathbb{C}^2$ 
such that $\omega_0(X,Y)>0$ is convex and hence contractible
and therefore locally near the boundary
this graph can be deformed to a graph which 
coincides with the annulus $A$ near $\gamma$ and hence
is a boundary
regular and boundary holomorphic immersed
symplectic disc with boundary $\gamma$ whose tangential
index coincides with the tangential index of $f$.
\end{proof}

\section{Boundaries of compact convex bodies with controlled curvature}

In this section we investigate periodic Reeb orbits on the 
boundary $\Sigma$ of a compact strictly convex body 
$C\subset \mathbb{C}^2$. Our main goal is the proof of
Theorem \ref{thm2} from the introduction.

We begin with observing that Corollary \ref{inverse} can
be applied to periodic 
Reeb orbits on boundaries of compact convex bodies.

\begin{lemma}\label{sympdisc}
Let $\gamma$ be a periodic Reeb orbit 
on $\Sigma$. Then there is a boundary regular symplectic immersion  
$f:(D,\partial D)\to 
(C,\gamma)$.
\end{lemma}
\begin{proof}
Let $\gamma$ be a periodic Reeb orbit on 
the boundary $\Sigma$ of a compact
strictly convex body  $C\subset \mathbb{C}^2$. Choose two
distinct points $a\not=b$ on $\gamma$ and smooth parametrizations
$\gamma_1,\gamma_2:[0,\pi]\to \gamma$ of the two subarcs of $\gamma$
connecting $a$ to $b$. We assume that the orientation
of $\gamma_2$ coincides with the orientation of $\gamma$ and that
the parametrizations $\gamma_1,\gamma_2$ coincide near $a,b$ with
the parametrization of $\gamma$ up to translation and 
reflection in the real line.
Define a map $f:(D,\partial D)\to (C,\gamma)$
as follows. Let $\tilde \gamma_1,\tilde \gamma_2:[0,\pi]\to
S^1$ be
parametrizations by arc length of the two half-circles 
of the unit circle $S^1\subset \mathbb{C}$ connecting
$1$ to $-1$, chosen in
such a way that
the orientation of $\tilde\gamma_2$ 
coincides with the orientation of $\partial D$. We
require that $f$ maps the line segment in $D$
connecting $\tilde \gamma_1(t)$ to $\tilde \gamma_2(t)$ which is
parametrized by arc length to the line segment in 
the convex body $C\subset \mathbb{C}^2$ connecting
$\gamma_1(t)$ to $\gamma_2(t)$ and parametrized propotional to arc
length on the same parameter interval. By construction, 
the map $f$ is smooth, moreover it is symplectic
near the points $1,-1$.

We claim that $f$ is a symplectic immersion.
For this let as before $\langle,\rangle$
be the usual euclidean inner product on $\mathbb{R}^4=\mathbb{C}^2$.
Let $t\in (0,\pi)$ and consider the
straight line segment $\ell$ in $C$ connecting $\gamma_1(t)$ to
$\gamma_2(t)$. By strict convexity of $C$, the arc $\ell$
is contained in $C$ and intersects $\Sigma$ transversely at the
endpoints. Let $X,Y$ be the tangents of $\ell$ at the endpoints 
$\gamma_1(t),\gamma_2(t)$ and
let as before $N$ be the outer normal field of $\Sigma$. 
Then
$\langle X,N(\gamma_1(t)) \rangle <0,\langle Y,N(\gamma_2(t)) \rangle >0$ 
and hence since
$\gamma_1^\prime(t)=-a_1JN(\gamma(t)),\gamma_2^\prime(t)
=a_2JN(\gamma_2(t))$ for some numbers $a_1>0,a_2>0$ 
we have 
$\omega_0(X,\gamma_1^\prime(t))>0$ and
$\omega_0(Y,\gamma_2^\prime(t))>0$. Now with respect to the 
usual trivialization of $T\mathbb{C}^2$ we have $X=Y$. On the
other hand, by the construction
of the map $f$, for every
point $s\in \ell$ the tangent space of $f(D)$ at $s$ is spanned by $X=Y$ 
and a convex linear combination of $\gamma_1^\prime(t),
\gamma_2^\prime(t)$. This shows that $f$ is a symplectic immersion.
Moreover $f$ is clearly boundary regular whence the lemma.
\end{proof}

We call an immersion $f:D\to C$ as in Lemma \ref{sympdisc}
a \emph{linear filling}
of the Reeb orbit $\gamma$. By Corollary \ref{inverse}, if 
${\rm lk}(\gamma)=-1$ then a linear filling $f$ of $\gamma$ satisfies
${\rm tan}(f)=0$. However, 
an immersed 
symplectic disc may have transverse self-intersection points
of negative intersection index, so 
there is no obvious relation between
the tangential index of a boundary
regular immersed symplectic disc and the number of
its self-intersection points. On the other
hand, if $\gamma$ admits an \emph{embedded} linear filling
then Corollary \ref{inverse} implies that 
${\rm lk}(\gamma)=-1$. 

Our final goal is to relate the Maslov index of a periodic 
Reeb orbit $\gamma$ to the geometry of the 
hypersurface $\Sigma$. 
For this consider for the moment an arbitrary bounded
domain $\Omega\subset \mathbb{C}^2$ with smooth 
boundary $\Sigma$ which is star-shaped with respect to the origin.
Write $C=\Omega\cup \Sigma$.
As before, denote
by $J$ the usual complex structure on $\mathbb{C}^2$ and 
let $\langle,\rangle$ be the \emph{euclidean} inner product. The
restriction $\lambda$ of the radial one-form $\lambda_0$ on
$\mathbb{C}^2$ defined by $(\lambda_0)_p(Y)= \frac{1}{2}\langle
Jp,Y\rangle$ $(p\in \mathbb{C}^2,Y\in T_p\mathbb{C}^2)$ defines a
smooth contact structure on $\Sigma$.

Let $N$ be the outer unit normal field of $\Sigma$.
As in Section 2 write 
$  \hat M(z_1,z_2)=(-\bar z_2,\bar z_1)$ and let 
$M$ be the section of $T\Sigma$ 
defined by $M(p)=\hat M\circ
N$. Its image is contained in 
the complex line subbundle $L$ of the tangent bundle
of $\Sigma$. The
sections $M,JM$ define a global trivialization of $L$ which is
symplectic with respect to the restriction of the symplectic form
$\omega_0$.

The kernel $\xi$ of the contact form is a smooth real 2-dimensional
subbundle of $T\Sigma$. Orthogonal projection $P$ of $T\Sigma$ onto
$L$ defines a smooth bundle epimorphism whose kernel is the
annihilator of the restriction of $\omega_0$ to $T\Sigma$. Thus 
the morphism $P$ preserves the restriction to $T\Sigma,L$ of the
symplectic form and therefore its restriction to the subbundle $\xi$
of $T\Sigma$ is a real symplectic bundle isomorphism.
Its inverse
$\pi:L\to \xi$ is a symplectic bundle morphism as well. Since by construction the
sections $M,JM$ of $L$ form a symplectic basis of $L$ we have

\begin{lem}\label{trivial}
The smooth sections $\pi\circ M,\pi\circ JM$ of the bundle $\xi$
define a symplectic trivialization $T:\xi\to (\mathbb{R}^2,dx\wedge
dy)$.
\end{lem}

In other words, for each $p\in \Sigma$ the restriction $T_p$ of $T$
to $\xi_p$ is an area preserving linear map $T_p:(\xi_p,\omega_0)\to
(\mathbb{R}^2,dx\wedge dy)$.

Recall from Section 2 that the Reeb vector field $X$ on $\Sigma$ is
given by
\begin{displaymath}
  X(p)=\phi(p)JN(p)
\end{displaymath}
where
\[\phi(p)=\frac{2}{\langle p,N(p) \rangle}>0.\]
Denote by $\Psi_t:\Sigma\to \Sigma$ the Reeb-flow of
$(\Sigma,\lambda)$ and let $\gamma$ be a periodic orbit for $\Psi_t$
of period $\chi >0$. 
Using the above trivialization $T$ of the bundle
$\xi$ we obtain a curve $\Phi:[0,T]\to SL(2,\mathbb{R})$ 
with $\Phi(0)={\rm Id}$ by defining
\begin{displaymath}
  \Phi(t):=T_{\Psi_t(p)}\circ d\Psi_t(p)\circ T_p^{-1}.
\end{displaymath}
where $p=\gamma(0)$. 
If the curve $\Phi$ is \emph{non-degenerate}, which means
that $\Phi(T)$ does not have one as an eigenvalue, then 
the Maslov index $\mu(\gamma)$ of $\gamma$ is defined as
the \emph{$\mu$-index} $\mu(\Phi)$ of the curve $\Phi$
as defined in \cite{HWZ95} .

To estimate the $\mu$-index of $\Phi$ define 
for a unit vector $X\in S^1\subset \mathbb{R}^2$ 
the \emph{rotation} of $X$ with respect to the curve $\Phi$
as the total rotation angle ${\rm rot}(\Phi,X)$ 
(or the total winding)
of the curve
\[t\to \frac{\Phi(t)X}{\Vert \Phi(t)X\Vert}\in S^1.\] 
The following lemma is valid for \emph{any} path
in $SL(2,\mathbb{R})$ beginning at the identity.
It uses an extension of the 
Maslov index to degenerate paths which is given in the proof
of the lemma.

\begin{lemma}\label{rotationindex}
Let $c:[0,\chi]\to SL(2,\mathbb{R})$ be a continuous arc
with $c(0)={\rm Id}$. Then
${\rm rot}(c,X)<(\mu(c)+1)\pi$ for every $X\in S^1$.
\end{lemma}
\begin{proof} We follow \cite{RS93}. 
Assume that $\mathbb{R}^2$ is equipped with the standard symplectic form.

In standard euclidean coordinates let $V=\mathbb{R}\times \{0\}\subset
\mathbb{R}^2$. The \emph{Maslov cycle} determined by $V$ is just
$V$, viewed as a point in the real projective line $\mathbb{R}P^1$ 
of all one-dimensional
subspaces of $\mathbb{R}^2$. A \emph{crossing} of a smooth curve   
$\Lambda:[a,b]\to \mathbb{R}P^1$ is a number
$t\in [a,b]$ such that $\Lambda(t)=V$. Then locally near
$t$, we can write $\Lambda(t)=\{x+A(s)x\}$ where
$A(s):V\to V^\perp$ is linear and vanishes for $s=t$ (and where
$V^\perp$ is the euclidean orthogonal complement of $V$ in $\mathbb{R}^2$).
With respect to the standard basis of $\mathbb{R}^2=V\oplus V^\perp$ 
we can view $s\to A(s)$  
as a real valued function. With this interpretation, 
the crossing is \emph{non-degenerate} if 
$A^\prime(t)\not=0$. The \emph{sign} ${\rm sign}\,\Gamma(\Lambda,V,t)$ 
of the crossing point $t$ then equals
the sign of $A^\prime(t)$ (p.830 of \cite{RS93}).
The \emph{Maslov index} of the curve $\Lambda:[a,b]\to \mathbb{R}P^1$ 
with only non-degenerate crossings is then defined to be
\[\mu(\Lambda,V)=
\frac{1}{2}{\rm sign}\,\Gamma(\Lambda,V,a)+
\sum_{a<t<b}{\rm sign}\,\Gamma(\Lambda,V,t)+
\frac{1}{2}{\rm sign}\,\Gamma(\Lambda,V,b)\]
(see p.831 of \cite{RS93}).

A smooth path $\Lambda:[a,b]\to \mathbb{R}P^1$
with $\Lambda(a)=V$ and only non-degenerate crossings lifts
to a smooth path $\tilde \Lambda:[a,b|\to S^1\subset \mathbb{C}=\mathbb{R}^2$
beginning at
$(1,0)$. Crossings of $\Lambda$ are precisely those points 
$t\in [a,b]$ where $\tilde \Lambda(t)=(\pm 1,0)$, and the sign
of the crossing is the sign of the derivative of $\tilde \Lambda$ 
with respect to the usual counter-clockwise orientation of $S^1$.
As a consequence, we have $\mu(\Lambda,V)=p$ for
$p\in \mathbb{Z}$ precisely if 
${\rm rot}(\tilde \Lambda)=p\pi$, and 
$\mu(\Lambda,V)=p+\frac{1}{2}$ precisely if 
${\rm rot}(\tilde \Lambda)\in (p\pi,(p+1)\pi)$. Here we write
${\rm rot}(\tilde \Lambda)$ to denote the total rotation of the path
$\tilde \Lambda$ in $S^1$.

For a curve $c:[0,\chi]\to SL(2,\mathbb{R})$, the Maslov index
$\mu(c,V)$ of $c$ with respect to $V$ is defined to be 
the Maslov index of the curve
$\Lambda:[0,\chi]\to \mathbb{R}P^1,t\to \Lambda(t)=
c(t)V$ \cite{RS93}
and hence ${\rm rot}(c,(1,0))< (\mu(c,V)+\frac{1}{2})\pi$.

The above definition of a Maslov index for paths 
in $SL(2,\mathbb{R})$ depends on the choice of the linear subspace $V$
(though this dependence can be removed by 
observing that any path in $SL(2,\mathbb{R})$ determines
a path of Lagrangian subspaces in $\mathbb{R}^4$, see \cite{RS93}).
To obtain an index for paths $c:[0,\chi]\to SL(2,\mathbb{R})$ 
beginning at $c(0)={\rm Id}$ 
which does not depend on such a choice
we proceed as follows.

A path $c:[0,\chi]\to SL(2,\mathbb{R})$ defines
a path $\alpha$ in the space of orientation preserving homeomorphisms
of $S^1$ by $\alpha(s)(X)=c(s)X/\Vert c(s)X\Vert$
$(s\in [0,\chi],X\in S^1)$. 
For each $s$ and each
$X$ we have $\alpha(s)(-X)=-\alpha(s)(X)$. 
This implies that 
$\vert {\rm rot}(c,X)-{\rm rot}(c,Y)\vert <\pi$ for any two
points $X,Y\in S^1$.
Now if there is some $p\in \mathbb{Z}$ such that
${\rm rot}(c,X)\in (2p\pi,2(p+1)\pi)$ for all $X\in S^1$ then 
we define $\mu(c)=2p+1$. By continuity, otherwise 
there is some $X\in S^1$ and some $p\in \mathbb{Z}$ 
with ${\rm rot}(c,X)=2p\pi$. By the above discussion,
the number $p$ is unique and we  
define $\mu(c)=2p$. With this definition of a Maslov index
for paths in $SL(2,\mathbb{R})$  beginning at the identity,
the statement of the lemma is obvious.

The fundamental group of $SL(2,\mathbb{R})$ with basepoint the
identity is infinitely
cyclic and generated by the loop $t\to e^{2\pi it}$ (viewed as a loop in 
$U(1)\subset SL(2,\mathbb{R})$). In particular, there is a 
natural group 
isomorphism $\rho:\pi_1(SL(2,\mathbb{R}))\to \mathbb{Z}$.
We claim that the Maslov index defined in the previous
paragraph has the following properties.
\begin{enumerate} 
\item $\mu(c)$ only depends on the homotopy class of $c$ with fixed
endpoints.
\item If $\alpha\in \pi_1(SL(2,\mathbb{R}))$ and $c:[0,1]\to SL(2,\mathbb{R})$
is any arc then \[\mu(\alpha \cdot c)=2\rho(\alpha)+\mu(c).\]
\item $\mu(c^{-1})=-\mu(c)$. 
\item For a constant invertible symmetric 
matrix $S$ with $\Vert S\Vert<2\pi$ and
for $c(t)=\exp tJS$ $(t\in [0,1])$ we have 
\[\mu(c)=\frac{1}{2}{\rm signature}(S).\]
\end{enumerate}

To see that these properties indeed hold true, note that 
the first statement is immediate from the definition.
By definition, the Maslov index of the standard rotation 
$\alpha:t\to e^{2\pi i t}\in U(1)$
equals $\mu(\alpha)=2=2\rho(\alpha)$ and hence
the second
statement follows from the definition and the first since 
the arc $\alpha \cdot c$ is homotopic with fixed endpoints to 
the concatentation of $c$ with (a representative of) $\alpha$.  
To see the fourth statement, observe that 
since $S$ is symmetric by assumption, we have 
$\Vert S\Vert <2\pi$ if and only if the absolute values of the eigenvalues
of $S$ are smaller than $2\pi$. 

Now by Theorem 3.2 of \cite{HWZ95},
the above properties uniquely determine the 
Maslov index of paths in $SL(2,\mathbb{R})$ beginning
at the identity and ending at a matrix which does not
have one as an eigenvalue as used
by Hofer, Wysocki and Zehnder (Section 3 of \cite{HWZ95},
see also Chapter 2 of \cite{S99} for alternative 
definitions). 
Together this completes the proof of the lemma.
\end{proof}

To calculate the total rotation angle of a vector under the curve 
$\Phi:[0,T]\to SL(2,\mathbb{R})$
defined above 
we use complex coordinates and 
view the vector fields $N,M$ as $\mathbb{C}^2$-valued functions on
$\Sigma$. For a curve $\gamma:[0,b]\to \Sigma$ we abbreviate
$N(t)=N(\gamma(t))$ and $M(t)=M(\gamma(t))$.  
Define a $U(2)$-valued
curve $O:[0,\chi]\to U(2)$ by the requirement that for each $t\in
[0,\chi]$, $O(t)$ is given with respect to the standard basis of
$\mathbb{C}^2$ by the matrix
\begin{displaymath}
  O(t):=\left(N(t),M(t)\right),\quad\mbox{ i.e. } O(t)
  \begin{pmatrix}
    a\\b
  \end{pmatrix}=aN(t)+bM(t) \quad \text{for}\, a,b\in\mathbb{C}.
\end{displaymath}
The image of the complex line $\{0\}\times \mathbb{C}\subset
\mathbb{C}^2$ under the map $O(t)$ is just the complex line
$L(\gamma(t))$. Therefore for each $t$, $\pi\circ O(t)$ is an
$\mathbb{R}$-linear isomorphism of $\{0\}\times \mathbb{C}$ onto
$\xi(\gamma(t))$.

From now on we use coordinates in $\IC^2$. Without loss of
generality we can assume that $\xi_{\gamma(0)}=L_{\gamma(0)}$ and
hence we have 
$\pi M(0)=M(0)=M$. 
%The Maslov index of $\Phi$ is then the
%variation of the argument of
%\[O(t)^{-1}\pi^{-1}d\Psi_t M_\alpha:[0,T]\to\IC-\{0\}.\]
If we define a unitary $2\times 2$ matrix $U(t)$ as
  \begin{equation}\label{u(t)}
    U(t):=\left(N(t),\frac{\pi^{-1}d\Psi_t M}
{\Vert \pi^{-1}d\Psi_t M\Vert}\right),
  \end{equation}
then the turning angle of $O(t)^{-1}\pi^{-1}d\Psi_t M$ about zero
(i.e. the rotation of $M$) 
is just the argument of $\det(U(t))$.

Define a unit vector field $\tilde M$ along $\gamma$ by
\[\tilde M(t)=\frac{\pi^{-1}d\Psi_tM}{\Vert
\pi^{-1}d\Psi_tM\Vert }.\] The following lemma is the main
technical tool for a calculation of the Maslov index of $\gamma$.
For its formulation, recall that 
the \emph{second fundamental form} of the hypersurface $\Sigma$ in
$\mathbb{C}^2$ is the symmetric bilinear form $\Pi:T\Sigma\times
T\Sigma\to \mathbb{R}$ which is defined as follows. Let $X,Y$ be
vector fields on $\Sigma\subset\IR^4$; then
\begin{displaymath}
  \Pi(X,Y)=-\langle dY(X),N \rangle=\langle dN(X),Y \rangle.
\end{displaymath}
The \emph{shape operator} of $\Sigma$ is the section $A$ of the
bundle $T^*\Sigma\otimes T\Sigma$ defined by $\Pi(X,Y)=\langle
AX,Y\rangle$.

\begin{lem}\label{derivative}
For $t_0\in [0,T]$ we have
\begin{align}
\frac{\partial}{\partial t}\det(U(t))|_{t=t_0}=
\frac{\partial}{\partial t}O(t)^{-1}\tilde M(t)\vert_{t=t_0} \notag\\
=i\phi(\gamma(t_0))(\Pi(JN(t_0),JN(t_0))+ \Pi(\tilde M(t_0),\tilde
M(t_0)))\det(U(t)) \notag
\end{align}
with $\phi(\gamma(t_0)) =\Vert X(\gamma(t_0))\Vert= \Vert
\gamma^\prime(t_0)\Vert=\frac{2}{\langle p,N(t_0) \rangle}$.
\end{lem}
\begin{proof}
In the sequel we always view the second fundamental form $\Pi$ of
$\Sigma$ as a bilinear form on a subspace of $\mathbb{R}^4$. Let
$\pi_2:\mathbb{C}^2\to \{0\}\times \mathbb{C}$ be the orthogonal
projection. Using the simple fact that
 \begin{displaymath}
   O^{-1}\circ\pi^{-1}=\pi_2\circ O^{-1}
 \end{displaymath}
we deduce
 \begin{align}\label{glei}
\frac{\partial}{\partial t} O(t)^{-1}\tilde M(t)& =
\frac{\partial}{\partial t} \frac{O(t)^{-1}\pi^{-1}d\Psi_t M}
{\Vert \pi^{-1}d\Psi_tM\Vert}\vert_{t=t_0}
\\ & =
\pi_2\bigl(\frac{\partial}{\partial t}O(t)^{-1}|_{t=t_0}\bigr)
\frac{d\Psi_{t_0} M}{\Vert \pi^{-1}d\Psi_{t_0} M\Vert}\notag \\
& +\pi_2 O(t_0)^{-1}\bigl(\frac{\partial}{\partial t}\left.
\frac{1}{\Vert\pi^{-1}d\Psi_t M\Vert}\right|_{t={t_0}}
\bigr)d\Psi_{t_0}M\notag \\
&+\pi_2 O(t_0)^{-1}\frac{1}{\Vert \pi^{-1}d\Psi_{t_0} M\Vert}
\bigl(\frac{\partial}{\partial t}d\Psi_{t}M|_{t={t_0}}\bigr).
\notag
\end{align}

The first term in our equation can be rewritten as
\begin{equation}
\begin{split}
 \pi_2\bigl(\frac{\partial}{\partial t} O(t)^{-1}|_{t=t_0}\bigr)&
\frac{d\Psi_{t_0} M}{\Vert\pi^{-1}d\Psi_{t_0} M\Vert}\\ 
= -\pi_2 O(t_0)^{-1} \bigl(\frac{\partial}{\partial t}
O(t)|_{t=t_0}O(t_0)^{-1} \bigr)&
\frac{d\Psi_{t_0} M}{\Vert \pi^{-1}d\Psi_{t_0} M\Vert}\notag \\
  &=*.
\end{split}
\end{equation}

By definition, for every $t$ the vectors $\{N(t),\tilde M(t)\}$ form
a unitary basis of $\IC^2$. Since $\{N(t),M(t)\}$ is also such a
unitary basis, there is a smooth function
$\psi:[0,\chi]\to\mathbb{R}$ such that $\tilde
M(t)=e^{i\psi(t)}M(t)$ for all $t$. We now use the second
fundamental form $\Pi$ to calculate the differential of the matrix
valued curve $O(t)=(N(t),M(t))$. For this recall the definition of
the shape operator $A:T\Sigma\to T\Sigma$ of $\Sigma$ and of the
orthogonal projection $P:T\Sigma\to L$. Since 
$M(t)=\hat M N(t)$ and
$\gamma^\prime(t)=\phi(t)JN(t)$ we have
\begin{equation}
\frac{\partial}{\partial t}O(t)=\phi(t)(AJN(t), 
\hat M AJN(t)) \notag.
\end{equation}
Thus with respect to the complex basis $(N(t),M(t))$ of
$\mathbb{C}^2$ we have
\begin{equation}
\frac{\partial}{\partial t}O(t)= \phi(t)\begin{pmatrix}
i\Pi(JN(t),JN(t)) & -\overline{PA(JN(t))}\notag\\
PA(JN(t)) & -i\Pi(JN(t),JN(t))
\end{pmatrix}
\end{equation}
where we view $PA(\phi JN(t))$ as a complex multiple of $M(t)$.

By the definition of the function $\psi$ we can write
\[\frac{O(t_0)^{-1}d\Psi_{t_0} M}{\Vert\pi^{-1}d\Psi_{t_0} M\Vert}=
\begin{pmatrix}
  ic(t_0)\\
  e^{i\psi(t_0)}
\end{pmatrix}\, \mbox{ for some } c(t_0)\in\IR.\]
Since $\tilde M(t)=e^{i\psi(t)} M(t)$ and
$O^{-1}(t)M(t)=
\left(\begin{smallmatrix}0\\1\end{smallmatrix}\right),
O^{-1}N(t)=\left(\begin{smallmatrix}1\\0\end{smallmatrix}\right)$ we
deduce that
\begin{align}
  *= & -\pi_2
\phi(t_0) \begin{pmatrix}
    -c(t_0)\Pi(JN(t_0),JN(t_0))-e^{i\psi(t_0)}\overline{PA(JN(t_0))} \notag\\
    ic(t_0)PA(JN(t_0))-ie^{i\psi(t_0)}\Pi(JN(t_0),JN(t_0))
  \end{pmatrix}.
\end{align}

Now let $\exp$ be the exponential map of the hypersurface $\Sigma$
with respect to the Riemannian metric induced by the euclidean
metric. We use this exponential map to compute the third term in
\eqref{glei}.
\begin{align}
\pi_2 O(t_0)^{-1} \frac{\frac{\partial}{\partial
t}d\Psi_{t}M |_{t={t_0}}} {\Vert\pi^{-1}d\Psi_{t_0} M\Vert}
=\frac{\pi_2 O(t_0)^{-1}} {\Vert\pi^{-1}d\Psi_{t_0} M\Vert}
\frac{\partial}{\partial t} \frac{\partial}{\partial
s}\Psi_t(\exp_{\gamma(0)}(sM))|_{t={t_0},s=0} \notag
\\
=\frac{\pi_2O(t_0)^{-1}} {\Vert\pi^{-1}d\Psi_{t_0} M\Vert}
\frac{\partial}{\partial s} (\phi(\Psi_{t_0}(\exp_{x_0}(sM))
JN(\Psi_{t_0}(\exp_{\gamma(0)}(sM))))\notag\\
=\pi_2O^{-1}(t_0) \left(\frac{1}{\Vert \pi^{-1}d\Psi_{t_0}M\Vert}
(\frac{\partial}{\partial s} \phi(\Psi_{t_0}(sM))JN(t_0)\vert_{s=0})
+\phi JA(d\Psi_{t_0}M)\right) \notag\\
=\pi_2O^{-1}(t_0)\phi(t_0)JA(\tilde M(t_0)+c(t_0)JN)\notag\\
=\phi(t_0)e^{i\psi(t_0)}\bigl(i\Pi(\tilde M,\tilde M)- \Pi(\tilde
M,J\tilde M)+ic\Pi(JN,\tilde M)- c\Pi(JN,J\tilde M)\bigr).\notag
\end{align}

The tangent vector of a curve $c$ in $\IC$ with constant norm always
has the form $\dot c=irc$ for some $r\in\IR$, so we can neglect the
radial parts of the above equations. Summing up the three terms in
\eqref{glei} yields
\begin{displaymath}
  \frac{\partial}{\partial t}
\det(U(t))|_{t=t_0}=i\phi(\Pi(JN,JN)+\Pi(\tilde M,\tilde M))
\det(U(t_0)).
\end{displaymath}
\end{proof}

The following corollary is immediate from Lemma \ref{derivative} and
the definition of the rotation of a vector with respect to an arc
in $SL(2,\mathbb{R})$.

\begin{corollary}\label{rotationpart}
Let $\gamma$ be a closed Reeb-orbit on $\Sigma$ with period $T$.
Then
\[{\rm rot}(\Phi,M(0))=
\int_0^T\vert  \gamma^\prime\vert
(\Pi(JN,JN)+\Pi(\tilde M,\tilde M))dt\]
where $\tilde M(t)=\frac{\pi^{-1}d\Psi_tM(0)}
{\Vert \pi^{-1}d\Psi_t M(0)\Vert}.$ 
\end{corollary}

Now we specialize again to the case that $C$ is a compact strictly convex
body in $\mathbb{C}^2$ with smooth boundary $\Sigma$ which
contains the origin in its interior. Recall that the \emph{total
curvature} of a smooth curve $\gamma:[0,t]\to \mathbb{C}^2$
parametrized by arc length is defined by
\[\kappa(\gamma)=\int_0^T\Vert \gamma^{\prime\prime}(t)\Vert dt.\]
The next corollary is immediate from Lemma \ref{rotationindex}
and lemma \ref{derivative}.

\begin{corollary}\label{totalcurvature}
Let $\Sigma$ be the boundary of a compact strictly convex body
$C\subset \mathbb{C}^2$. If the principal curvatures $a\geq b\geq c$
of $\Sigma$ satisfy the pointwise 
pinching condition $a\leq b+c$ then the
total curvature of a periodic Reeb orbit of Maslov index 3 is
smaller than $4\pi$.
\end{corollary}
\begin{proof}
Let $\Sigma$ be the boundary of a compact 
strictly convex body $C\subset
\mathbb{C}^2$ with principal curvatures $a\geq b\geq c$ satisfying
the pinching condition $a\leq b+ c$. Let $\gamma:[0,T]\to \Sigma$
be a periodic Reeb orbit of Maslov index 3. We assume that $\gamma$
is parametrized by arc length on $[0,T]$. By 
Lemma \ref{rotationindex}, the rotation of the vector 
$M(\gamma(0))$ under the derivative of the Reeb flow
is smaller than $4\pi$.

Denote by $N(t)$ the normal field of the sphere restricted to the
curve $\gamma$. Then $\gamma^\prime(t)=JN(t)$ and therefore
\[\kappa(\gamma)=\int_0^T
\Vert \frac{\partial}{\partial t}N(t)\Vert dt \leq \int_0^T
a(\gamma(t))dt \leq \int_0^Tb(\gamma(t))+ c(\gamma(t))dt<4\pi\] by
Corollary \ref{rotationpart}.
\end{proof}

We use Corollary \ref{totalcurvature} to complete the proof of 
Theorem \ref{thm2} from the introduction.

\begin{proposition}\label{selflink}
Let $\Sigma$ be the boundary of a compact strictly convex domain
$C\subset \mathbb{C}^2$. If the principal curvatures $a\geq b\geq c$
of $\Sigma$ satisfy the pointwise 
pinching condition $a\leq b+c$ then a
periodic Reeb orbit on $\Sigma$ of Maslov index 3 
bounds an embedded symplectic disc $f:(D,\partial D)\to 
(C,\gamma)$. In particular, $\gamma$ 
has self-linking
number $-1$.
\end{proposition}
\begin{proof}
Define the \emph{crookedness} of a smooth closed curve
$\gamma:S^1\to \mathbb{R}^4$ to be the minimum of the numbers
$m(\gamma,v)$ where $m(\gamma,v)$ is the number of minima of the
function $t\to \langle\gamma(t),v\rangle,v\in S^3$. By a result of Milnor
\cite{milnor}, the crookedness of a curve of total curvature smaller
than $4\pi$
equals one. Thus by Corollary \ref{totalcurvature}, if
$\gamma$ is a periodic Reeb orbit on $\Sigma$ of Maslov index 3 then 
there is some $v\in
S^3$ such that the restriction to $\gamma$ 
of the function $\phi:x\to \langle x,v\rangle$
assumes precisely one maximum and one minimum.
We may moreover assume that these are the only critical points
of the restriction of $\phi$ to $\gamma$ and that they
are non-degenerate (see \cite{milnor}). 

Let $a$ be the unique minimum of $\phi$ 
on $\gamma$. Assume that $\gamma$ is parametrized
in such a way that
$\gamma(0)=a$. Let $\gamma_2:[0,\sigma]\to \Sigma$
be the parametrized subarc of $\gamma$ issuing from $\gamma_2(0)=a$
which connects $a$ to the unique 
maximum $b$ of $\phi$ on $\gamma$. Let $\gamma_1:[0,\sigma]\to \Sigma$
be the parametrization of the second subarc of $\gamma$ connecting
$a$ to $b$ such that $\phi(\gamma_1(t))=\phi(\gamma_2(t))$ 
for all $t\in [0,\sigma]$; 
this is possible by construction
and by our choice of $\phi$. Then the symplectic disc obtained from
this parametrization by linear
filling as in the proof of Lemma 4.1
is embedded. By Corollary \ref{inverse},
this implies that the self-linking number of
$\gamma$ equals $-1$.
\end{proof}

{\bf Acknowledgement:} The authors thank the referee of an earlier
version of this paper for pointing out a gap in the proof of
Theorem 2 and for suggesting the statement of Theorem 1.

\noindent
MATHEMATISCHES INSTITUT DER UNIVERSIT\"AT\\
ENDENICHER ALLEE 60\\
53115 BONN GERMANY

\end{document}